\documentclass{siamltex}
\usepackage[dvips]{graphicx}
\usepackage{amsmath}
\usepackage{latexsym}
\usepackage{amssymb}
\usepackage{color}
\usepackage{algorithm,algorithmic}
\usepackage{epsfig}
\usepackage{multirow}
\usepackage{graphicx}
\usepackage{caption}
\usepackage{subcaption}

\newtheorem{remark}{Remark}[section]

\makeatletter
\newcommand{\ssymbol}[1]{^{\@fnsymbol{#1}}}
\makeatother

\newcommand{\R}{\mathbb R}

\usepackage[mathscr]{euscript}

\title{Discrete cosine transform LSQR and GMRES methods for  multidimensional ill-posed problems}

\author{A. El Ichi\footnotemark [1]\thanks{Laboratoire de Mathématiques, Informatique et Applications, S\'ecurit\'e de l'Information LABMIA-SI, University Mohamed V, Rabat Morocco }  \and  M. El Guide \thanks{Centre for Behavioral Economics and Decision Making(CBED), FGSES, Mohammed VI Polytechnic University, Green City, Morocco} \and   K. Jbilou\footnotemark[1] \thanks{LMPA, 50 rue F. Buisson, ULCO Calais, France; Mohammed VI Polytechnic University, Green City, Morocco; jbilou@univ-littoral.fr }}

\date{}

\begin{document}
	
	%%%%%%%%%%%%%%%%%%%%%%%%%
	
	%
	% % !TeX spellcheck = en_GB
	
	%
	%
	%
	% %--------numbring equations-----------
	% %\numberwithin{equation}{section}
	% %-------------------------------------
	%
	% %%%%%%%%%%%%%%%%%%%%%%%%%%%%%%%%%%%%%%%%%%%%%
	% % \newtheorem{proposition}{Proposition}[section]
	% % \newtheorem{theorem}{Theorem}[section]
	% % \newtheorem{definition}{Definition}[section]
	% % \newtheorem{lemma}{Lemma}[section]
	% % \newtheorem{corollary}{Corollary}[section]
	% %\newtheorem{acknowledgement}[theorem]{Acknowledgement}
	%% \newtheorem{remark}{Remark}[section]
	%
	% %%%%%%%%%%%%%%%%%%%%%%%%%%%%%%%%%%%%%%%%%%%%%%%%%%%%%%%%
	% %\parindent 0.5cm
	%\textwidth 17cm
	% %\textheight 20cm
	% %\evensidemargin -0.1cm
	% %\oddsidemargin -0.1cm
	%

	\maketitle

	%	\date{Received: date / Accepted: date}
	% The correct dates will be entered by the editor

	%\maketitle
	
	\begin{abstract}    
		In the present work, we propose new tensor Krylov subspace method for ill posed linear tensor problems such as in color or video image restoration. Those methods are based on the tensor-tensor discrete cosine transform that gives fast tensor-tensor product computations. In particular, we will focus on the tensor discrete cosine versions of GMRES, Golub-Kahan bidiagonalisation and LSQR methods. The presented numerical tests show that the methods are very fast and give good accuracies when solving some linear tensor ill-posed problems. 	
	\end{abstract}


	% \end{frontmatter}

	\medskip
	
	%====================================================================
	
	%===========
	%\vspace*{1cm}*

	\section{Introduction}
	The aim of this paper is to solve the following tensor problem

	\begin{equation}\label{eq1}
\displaystyle \min_{\mathscr{X}} \Vert 	{\mathcal M} (\mathscr{X})-  \mathscr{C} \Vert_F
	\end{equation}
	where  ${\mathcal M}$ is a linear operator  that could be described as 
	\begin{equation}\label{eq2}
	{\mathcal M} (\mathscr{X}) = \mathscr{A} \star_c\mathscr{X},\; {\rm or}\; 	{\mathcal M} (\mathscr{X}) = \mathscr{A} \star_c\mathscr{X}\star_c\mathscr{B},
	\end{equation}
where $\mathscr{A} \in \mathbb{R}^{n_1 \times n_2 \times n_3}$ is a three mode tensor, $\mathscr{X} \in \mathbb{R}^{n_2 \times s \times n_3}  $, $\mathscr{B}\in \mathbb{R}^{s \times s \times n_3} $ and $\mathscr{C}   \in \mathbb{R}^{n_1 \times s \times n_3}$ are  three mode tensors,   and $\star_c$ is the cosine product to be also defined later. Applications of such  problems arise in  signal processing \cite{lb},  data mining \cite{lxnm},  computer vision and so many other modern applications  in machine learning.  For large scale problems,  we have to  take advantage of the  multidimensional structure to build rapid and robust iterative methods. Tensor Krylov subspace methods could be useful and very fast solvers for those tensor problems. \\ In the present paper,  we will be  interested
	in  developing robust and fast iterative  tensor  Krylov based subspace methods using  tensor-tensor products such as the tensor cosine product \cite{ahmad}.  In many applications such as in image or video processing, the obtained discrete problems are very ill conditioned and the we should add some regularization techniques such as the generalized cross validation method. Standard and global Krylov subspace methods are suitable  when dealing with grayscale images, e.g, \cite{belguide,reichel2,reichel1}. However, these methods might be time consuming to numerically solve problems related to multi channel images (e.g. color images, hyper-spectral images and videos). 
	
%	For the Einstein product, both the Einstein tensor global   Arnoldi and Einstein
%	tensor global  Gloub-Kahan bidiagonalization algorithms have been established \cite{Elguide}, which makes so  natural to generalize these methods using the T-product. 
	 In this paper, we will show that the
	tensor-tensor product  between third-order
	tensors allows the application of the global iterative methods, such as the global Arnoldi and global Golub-Kahan algorithms. The tensor form of the proposed Krylov methods, together with using the fast
cosine  transform (DCT) to compute the c-product between third-order
	tensors can be efficiently implemented on many modern computers and allows to significantly reduce the overall computational complexity. It is also worth mentioning that our approaches can be naturally generalized to higher-order tensors in a recursive manner.
	
\noindent  This paper is organized as follows. We shall first present in Section 2 some symbols and notations used throughout the paper. We also recall   some definitions related to the cosine product between two tensors.   In Section 3, we present some inexpensive approaches
based on cosine global Krylov subspace methods combined with regularization techniques to
solve the obtained ill-posed  tensor problem (\ref{eq1}) . Section 5 is dedicated to some numerical experiments.% Concluding
	\section{Definitions and Notations} A tensor is  a multidimensional array of data. The number of indices of a tensor is called modes or ways. 
	Notice that a scalar can be regarded as a zero mode tensor, first mode tensors are vectors and matrices are second mode tensor. The order of a tensor is the dimensionality of the array needed to represent it, also known as
	ways or modes. 
	For a given N-mode (or order-N) tensor $ \mathscr {X}\in \mathbb{R}^{n_{1}\times n_{2}\times n_{3}\ldots \times n_{N}}$, the notation $x_{i_{1},\ldots,i_{N}}$ (with $1\leq i_{j}\leq n_{j}$ and $ j=1,\ldots N $) stands for the element $\left(i_{1},\ldots,i_{N} \right) $ of the tensor $\mathscr {X}$.  \\  Fibers are the higher-order analogue of matrix rows and columns. A fiber is
	defined by fixing all the indexes  except  one. A matrix column is a mode-1 fiber and a matrix row is a mode-2 fiber. Third-order tensors have column, row and tube fibers. An element $c\in \mathbb{R}^{1\times 1 \times n}$ is called a tubal-scalar of length $n$. More details are found in  \cite{kimler1,kolda1}.  \\ In the present paper, we will consider only 3-order tensors and show how to use them in color image and video processing.

	\subsection{Discrete Cosine Transformation}
 \noindent In this subsection we recall some definitions and properties of the discrete cosine transformation and the c-product. The Discrete Cosine Transformation (DCT) plays a very important role in the definition of the c-product of tensors. The DCT on 
 	a vector $v \in {\R}^n$ is defined by
 	\begin{equation}
 	\label{dft1}
 	\tilde v= C_n v \in {\R}^n,
 	\end{equation}
 	where $C_n$ is the $n\times n$  discrete cosine
 	transform matrix with entries 
 	$$ (C_n)_{ij}=\sqrt{\frac{ 2-\delta_{i1}}{n}}\cos\left( \dfrac{(i-1)(2j-1)\pi}{2n}\right)\quad 1<i,j<n $$
 	with $\delta_{ij} $ is the Kronecker delta for more details. Its  known that the matrix $C_n$ is orthogonal, i.e, $C_n^TC_n=C_nC_n^T=I_n$; see \cite{Ng}. Furthermore, for  any  vector $v \in \mathbb{R}^n$, the   matrix vector multiplication $C_nv$     can be computed in $O(nlog(n))$ operations.  Also, Ng  and  al. \cite{Ng}  \textcolor{blue}{showed} that   matrices which  can be diagonalized by  $C_n$ are some
 	special Toeplitz-plus-Hankel matrices. In other words, we have 
 	 \begin{equation}\label{dft5}
 	C_n \, {\tt th }(v)\, C_n^{-1}  = {\rm Diag}(\tilde v),\\
 	\end{equation}
 %	which is equivalent to
 %	\begin{equation}\label{dft6}
 %	F_n \, {\rm circ}(v)\, F_n^{*}  = n\,{\rm Diag}(\tilde v),\\
 %	\end{equation}
 	where
 	$$
 	{\tt th}(v)=  \mathop{ {\underbrace{\left ( 
 				\begin{array}{cccc}
 				v_1 & v_2 & \ldots & v_n\\
 				v_2 & v_1 & \ldots & v_3\\
 				\vdots & \vdots  & \ldots & \vdots \\
 				v_n & v_{n-1} & \ldots & v_1\\
 				\end{array}
 				\right )}}}\limits_{\tiny{\text{Toeplitz}}}  +  \mathop{ {\underbrace{\left( {\begin{array}{*{20}{c}} {{v_2}}&\ldots &  {{v_n}} &0\\
 					\vdots &  \reflectbox{ $\ddots$}   & \reflectbox{ $\ddots$} &   {{v_{{n}}}}\\
 					{{v_{n}}}&0& \ldots&  \vdots\\
 					0&{{v_{{n}}}}&\ldots&  {{v_2}}
 					\end{array}} \right)}}}\limits_{\tiny{\text{Hankel}}}   
 	$$
 	and $ {\rm Diag}(\tilde v)$ is the diagonal matrix whose $i$-th diagonal element is $ (\tilde v)_i$.

	\subsection{Properties of the cosine product}
    	
\medskip
\noindent In this subsection, we briefly review some concepts and notations,  that play a central role for the elaboration of the tensor iterative methods based on the c-product; see \cite{klimer5} for more details on the c-product.

	Let $\mathscr {A} \in \mathbb{R}^{n_{1}\times n_{2}\times n_{3}} $ be a third-order tensor, then the operations {\tt mat}  and  its reverse {\tt ten} are defined by
	\begin{align*}
 	  {\tt mat}(\mathscr {A})=   \mathop{ {\underbrace{\left ( 
				\begin{array}{cccc}
				A_1 & A_2 & \ldots & A_n\\
				A_2 & A_1 & \ldots & A_3\\
				\vdots & \vdots  & \ldots & \vdots \\
				A_n & A_{n-1} & \ldots & A_1\\
				\end{array}
				\right )}}}\limits_{\tiny{\text{Block Toeplitz}}}  &+  \mathop{ {\underbrace{\left( {\begin{array}{*{20}{c}} {{A_2}}&\ldots &  {{A_n}} &0\\
					\vdots &  \reflectbox{ $\ddots$}   & \reflectbox{ $\ddots$} &   {{A_{{n}}}}\\
					{{A_{n}}}&0& \ldots&  \vdots\\
					0&{{A_{{n}}}}&\ldots&  {{A_2}}
					\end{array}} \right)}}}\limits_{\tiny{\text{Block  Hankel}}}  \in {\R}^{ n_1n_3 \times n_2n_3}
				\end{align*}				
			and the reverse operation denoted by  {\tt ten} and such that 
   $$ {\tt ten}({\tt mat}(\mathscr {A}) ) =  \mathscr {A}.$$
	\noindent Let $\widetilde {\mathscr{A}}$ be the tensor obtained by applying the DCT on all the tubes of the tensor $\mathscr {A}$. With the Matlab command ${\tt dct}$ as 
	$$\widetilde {\mathscr{A}}= {\tt dct}(\mathscr {A},[ \;],3), \; {\rm and }\;\; {\tt idct} (\widetilde {\mathscr{A}}, [ \;],3)= \mathscr {A},$$
 	where ${\tt idct}$ denotes the Inverse Discrete Cosine Transform.\\
\begin{remark}
	 Notice that   the  tensor $\widetilde {\mathscr{A}}$ can be computed   by  using  the  3-mode product defined in   \cite{kolda1} as  follows:
	 $$ \widetilde{{\mathscr{A} }}= {\mathscr{A} }\times_3 M,    $$ 
	 where  M is the  ${ n_{3}\times n_{3}}$  invertible matrix  given by  $$M=W^{-1}C_{n_3}(I+Z),$$ and  $C_{n_3}$ denotes the $ n_3\times n_3 $ Discrete Cosine   Transform  DCT matrix, $W= {\rm diag}(C_{n_3}(:,1))$ is the diagonal matrix made of the first column of the DCT matrix, Z is an  $n_3\times n_3 $ circulant  matrix which   can be computed in MATLAB using the command $\tt W={\rm diag}({\rm ones}(n_{3}-1,1),1)$ and $I$ the $ n_3\times n_3 $ identity matrix; see \cite{klimer5} for more details.\\
\end{remark}  
 Let ${\bf A}$ be the matrix 
\begin{equation}\label{dft9}
{\bf A}= \left (
\begin{array}{cccc}
{A}^{(1)}& &&\\
& {A}^{(2)}&&\\
&&\ddots&\\
&&&{A}^{(n_3)}\\
\end{array}
\right)\in \mathbb{R}^{n_3 n_1\times n_3n_2 },
\end{equation}
where  the matrices ${A}^{(i)}$'s are the frontal slices of the tensor ${\widetilde {\mathscr{A}}}$. The  block   matrix ${\tt mat}(\mathscr {A})$ can also be block diagonalized using the DCT matrix and this gives
\begin{equation}\label{dft8}
(C_{n_3} \otimes I_{n_1})\, {\tt mat}(\mathscr {A})\, 	(C_{n_3}^{T} \otimes I_{n_2})={\bf A}.
\end{equation}
\begin{definition}
	The \textbf{c-product}  between  two tensors
	$\mathscr {A} \in \mathbb{R}^{n_{1}\times n_{2}\times n_{3}} $ and $\mathscr {B} \in \mathbb{R}^{n_{2}\times m\times n_{3}} $ is an  ${n_{1}\times m\times n_{3}}$ tensor  given by:	
	$$\mathscr {A} \star_c \mathscr {B}={\tt ten}({\tt mat}(\mathscr {A}){\tt mat}(\mathscr {B}) ).$$
\end{definition}
Notice that from the relation \eqref{dft9}, we can show that the   product $\mathscr {C}= \mathscr {A} \star_c \mathscr {B}$ is equivalent to $  {\bf C}= {\bf A}\,{\bf B}$. %So, the efficient way to compute the T-product is to use Discret Fourier Transform (DCT). \\
%Using the relation \eqref{f1}, 
The following algorithm allows us to compute, in an efficient way, the c-product of the tensors $\mathscr {A}$ and 
$\mathscr {B}$, see \cite{klimer5}.\\
\begin{algorithm}[!h]
	\caption{Computing the  c-product }\label{algo1}
\textbf{	Inputs}: $\mathscr {A} \in \mathbb{R}^{n_{1}\times n_{2}\times n_{3}} $ and $\mathscr {B} \in \mathbb{R}^{n_{2}\times m\times n_{3}} $\\
\textbf{	Output}: $\mathscr {C}= \mathscr {A} \star_c \mathscr {B}  \in \mathbb{R}^{n_{1}\times m \times n_{3}} $
	\begin{enumerate}
		\item Compute $ \widetilde{{\mathscr{A} }}= {\tt dct}(\mathscr {A},[ \;],3)$ and $\mathscr {\widetilde B}={\tt dct}(\mathscr {B},[\; ],3)$.
		\item Compute each frontal slices of $\mathscr {\widetilde C}$ by
		$$C^{(i)}= 
		A^{(i)} B^{(i)}  $$
		\item Compute $  {{\mathscr{C} }}= {\tt idct}(\mathscr {\widetilde C},[ \;],3)$  .	 	
	\end{enumerate}
\end{algorithm}
	\noindent For the c-product, we have the following definitions and remarks
\begin{definition} 
	The identity tensor $\mathscr{I}_{n_{1}n_{1}n_{3}} $ is the tensor such  that  all frontal slice  of   $\widetilde{{\mathscr{I} }}_{n_{1}n_{1}n_{3}}$  is the identity matrix $I_{n_1n_1}$ .\\
			An $n_{1}\times n_{1} \times n_{3}$ tensor $\mathscr{A}$ is invertible, if there exists a tensor $\mathscr{B}$ of order  $n_{1}\times n_{1} \times n_{3}$  such that
		$$\mathscr{A}  \star_c \mathscr{B}=\mathscr{I}_{ n_{1}  n_{1}  n_{3}} \qquad \text{and}\qquad \mathscr{B}  \star_c \mathscr{A}=\mathscr{I}_{ n_{1}  n_{1}  n_{3}}.$$
		In that case, we set $\mathscr{B}=\mathscr{A}^{-1}$. 	It is clear that 	$\mathscr{A}$ is invertible if and only if   ${\tt mat }(\mathscr{A})$ is invertible.\\
	The  inner scalar product is defined by
		$$\langle \mathscr{A}, \mathscr{B} \rangle = \displaystyle \sum_{i_1=1}^{n_1} \sum_{i_2=1}^{n_2}  \sum_{i_3=1}^{n_3} a_{i_1 i_2 i_3}b_{i_1 i_2 i_3}$$
and he  corresponding norm  is given by
		$ \Vert \mathscr{A} \Vert_F=\displaystyle \sqrt{\langle  \mathscr{A} ,  \mathscr{A}  \rangle}.$\\
			An $n_{1}\times n_{1} \times n_{3}$ tensor  $\mathscr{Q}$  is orthogonal if
		$\mathscr{Q}^{T}   \star_c  \mathscr{Q}=\mathscr{Q} \star_c \mathscr{Q}^{T}=\mathscr{I}_{ n_{1}  n_{1}  n_{3}}.$\\
\end{definition}
\medskip
\begin{remark}
	Another interesting way for computing the scalar product and the associated norm is as follows:
	$\langle \mathscr{A}, \mathscr{B} \rangle = \displaystyle \frac{1}{n_3}  \langle  {\bf A}, {\bf B} \rangle$ and $  \Vert \mathscr{A} \Vert_F= \displaystyle \frac{1}{\sqrt{n_3}} \Vert {{\bf A}} \Vert_F,$
	where the block diagonal matrix ${\bf A}$  is defined by \eqref{dft9}.
\end{remark}

\medskip
\noindent We now introduce the new c-diamond tensor-tensor product.
\begin{definition}%\cite{ElIchi}
		Let
		$\mathscr{A} =[\mathscr{A}_{1},\ldots,\mathscr{A}_{p}]\in {\mathbb R}^{n_{1}\times ps \times n_3},$ where $  \mathscr{A}_{i} \in {\mathbb R}^{n_{1}\times s \times n_3} , \, i =1,...,p$ %\hspace{6cm}
		and let $	\mathscr{B} =[\mathscr{B}_{1},\ldots,\mathscr{B}_{l}]\in {\mathbb R}^{n_{1}\times \ell s \times n_3}$ with $  \mathscr{B}_{j}\in {\mathbb R}^{n_{1}\times s \times n_3}, \, j =1,...\ell$. Then,  the product $\mathscr{A}^{T} \diamondsuit  \mathcal{B} $ is the   $p \times  \ell  $  matrix  given by : 
		$$ (\mathscr{A}^{T} \diamondsuit  \mathcal{B})_{i,j}  =   \langle \mathscr {A}_i  ,\mathscr{B}_{j} \rangle  \;\;.   
		$$ 
\end{definition}

	\section{Tensor discrete cosine global Krylov subspace methods}
	
	In this section, we propose iterative
	methods based on  tensor cosine global Arnoldi and  cosine global Golub– Kahan bidiagonlization (cosine-GGKB),
	combined with Tikhonov regularization, to solve  some discrete ill posed problems. We consider the following discrete ill-posed tensor equation
	\begin{equation}\label{tr1}
	\mathscr{A} \star_c \mathscr{X}= \mathscr{C},\quad \mathscr{C}=\widehat{ \mathscr C}+  \mathscr{N},
	\end{equation}
	where $\mathscr{A} \in {\R}^{n \times m \times p}$, $\mathscr{X}$,  $ \mathscr{N}$ (additive noise) and $\mathscr{C}$ are tensors in ${\R}^{n \times s \times p}$. 
	In color image processing, $p=3$, $\mathscr{A}$ represents the blurring tensor, $\mathscr{C}$ the blurry and noisy observed  image,  $\mathscr{X}$ is the image that we would like to restore  and  $\mathscr{N}$ is an unknown additive noise. Therefore, to stabilize  the recovered image, regularization techniques are  needed. There are several techniques to regularize the linear inverse problem given by
	equation (\ref{tr1}); for the matrix case, see for example, \cite{belguide,reichel2,golubwahba,hansen1}. All of these techniques stabilize the restoration process by adding a regularization term, depending on some  priori knowledge of  the unknown image. One of the most regularization method is due to Tikhonov and is given as follows 
	\begin{equation}\label{tr2}
	\underset{\mathscr{X}}{\text{min}}\{\|\mathscr{A} \star_c \mathscr{X} - \mathscr{C}  \|_F^2+\lambda \|\mathscr{X}\|_F^2\}.
	\end{equation}
	
	\noindent Many techniques for choosing a suitable value of $\lambda$  have been analysed and illustrated
	in the literature; see, e.g., \cite{CR,golubwahba,hansen1,wahbagolub}  and references therein. In this paper we will use the
	discrepancy principle and the Generalized Cross Validation (GCV) techniques.

\subsection{The tensor discrete cosine  GMRES }
\noindent Let $\mathscr{A}\in \mathbb{R}^{n\times n \times p}$ and $\mathscr{V}\in \mathbb{R}^{n\times s \times p}$. We introduce the  tensor Krylov subspace $\mathcal{\mathscr{TK}}_m(\mathscr{A},\mathscr{V} )$ associated to the cosine-product, defined for  the pair $(\mathscr{A},\mathscr{V})$   as follows
\begin{equation}
\label{tr3}
\mathcal{\mathscr{TK}}_m(\mathscr{A},\mathscr{V} )= {\rm Tspan}\{ \mathscr{V}, \mathscr{A} \star_c\mathscr{V},\ldots,\mathscr{A}^{m-1}\star_c\mathscr{V} \}\\
=\left\lbrace \mathscr{Z} \in \mathbb{R}^{n\times s \times p}, \mathscr{Z}= \sum_{i=1}^m \alpha_{i} \left(
\mathscr{A}^{i-1}\star_c\mathscr{V}\right) \right\rbrace,
\end{equation}
where $\alpha_{i}\in \mathbb{R} $,  $\mathscr{A}^{i-1}\star_c\mathscr{V}=\mathscr{A}^{i-2}\star_c\mathscr{A}\star_c\mathscr{V}$, for $i=2,\ldots,m$ and $\mathscr{A}^{0}$ is the identity tensor. In the following algorithm, we define the  Tensor cosine-global Arnoldi algorithm.

\begin{algorithm}[H]
	\caption{Tensor discrete cosine Arnoldi} \label{TGA}
	\begin{enumerate}
		\item 	{\bf Input.} $\mathscr{A}\in \mathbb{R}^{n\times n \times p}$, $\mathscr{V}\in \mathbb{R}^{n\times s \times p}$ and the positive integer $m$.
		\item Set $\beta=\|\mathscr{V}\|_F$, $\mathscr{V}_{1} =   \dfrac{\mathscr{V}}{  \beta}$
		\item For $j=1,\ldots,m$
		\begin{enumerate}
			\item $\mathscr{W}=  \mathscr{A}\star_c   \mathscr{V}_j$
			\item for $i=1,\ldots,j$
			\begin{enumerate}
				\item $h_{i,j}=\langle \mathscr{V}_i, \mathscr{W} \rangle$
				\item $\mathscr{W}=\mathscr{W}-h_{i,j}\;\mathscr{V}_i$	
			\end{enumerate}	
			\item End for
			\item   $h_{j+1,j}=\Vert \mathscr {W} \Vert_F$. If $h_{j+1,j}=0$, stop; else
			\item $\mathcal {V}_{j+1}=\mathscr {W}/h_{j+1,j}$.

		\end{enumerate}
		\item End
	\end{enumerate}
\end{algorithm}

\medskip
It is not difficult to show that after  $m$ steps of Algorithm \ref{TGA},  the tensors  $\mathscr{V}_{1},\ldots,\mathscr{V}_{m}$ form an orthonormal basis of the tensor   Krylov  subspace $\mathscr{TK} _{m}(\mathscr{A},\mathscr{V})$.
%	\end{proposition}
\noindent Let $\mathbb{V}_{m}  $ be  the $(n\times (sm)\times p)$ tensor with frontal slices $\mathscr{V}_{1},\ldots,\mathscr{V}_{m}$ and let $ {\widetilde{H}}_{m}$ be  the $(m+1)\times m  $ upper  Hesenberg matrix whose elements are the $h_{i,j}$'s defined by Algorithm \ref{TGA}.    Let  $ {H}_{m}$ be the matrix obtained from $\widetilde{ { H}}_{m}$ by deleting its last row;  $H_{.,j}$ will denote the $j$-th column of the matrix  $H_m$ and $\mathscr{A}\star_c\mathbb{V}_{m}  $ is  the $(n\times (sm)\times p)$ tensor with frontal slices  $\mathscr{A}\star_c\mathscr{V}_{1},\ldots,\mathscr{A}\star_c\mathscr{V}_{m}$: 
\begin{equation}
\label{ev1}
\mathbb{V}_{m}:=\left[  \mathscr{V}_{1},\ldots,\mathscr{V}_{m}\right] \;\;\; {\rm and}\;\;\; \mathscr{A}\star_c \mathbb{V}_{m}:=[\mathscr{A}\star_c\mathscr{V}_{1},\ldots,\mathscr{A}\star_c\mathscr{V}_{m}].
\end{equation}
We introduce the product  $\circledast$ defined by $$\mathbb{V}_{m}\circledast y=\sum_{j=1}^{m} {y}_{j}\mathscr{  {V}}_{j},\; \; y= (y_1,\ldots,y_m)^T\in \mathbb{R}^m,\, and \; 
\mathbb{V}_{m}\circledast {{ {H}}_{m}}=\left[   \mathbb{V}_m\circledast H_{.,1} ,\ldots,\mathscr{V}_{m}\circledast H_{.,m} \right].$$

\noindent With the above notations, we can easily prove the  results of the following proposition.

\begin{proposition}\label{T-GlobalArnolproposit}
	Suppose that m steps of Algorithm \ref{TGA} have been   run. Then, the following statements hold:
	\begin{eqnarray}
	\mathscr{A}\star_c\mathbb{V}_{m}&=&\mathbb{V}_{m}\circledast {{ {H}}_{m}} +  h_{m+1,m}\left[  \mathscr{O}_{n\times s\times p},\ldots,\mathscr{O}_{n\times s\times p},\mathscr{V}_{m+1}\right],\\
	\mathscr{A}\star_c\mathbb{V}_{m}&=&\mathbb{V}_{m+1}   \circledast \widetilde{ { H}}_{m}, \\
	\mathbb{V}_{m}^{T}\diamondsuit\mathscr{A}\star_c\mathbb{V}_{m}&=& {H}_{m}, \\	
	\mathbb{V}_{m+1}^{T}\diamondsuit  \mathscr{A}\star_c\mathbb{V}_{m}&=&\widetilde{ { H}}_{m},\\
	\mathbb{V}_{m}^{T} \diamondsuit\mathbb{V}_m&=& {I}_{ m  },\\
		\Vert  \mathbb {V}_{m}  \circledast y \Vert_F & = & \Vert y \Vert_2, \; y \in \mathbb{R}^m,
	\end{eqnarray}
	where ${I}_{ m  }$ the identity matrix and $\mathscr{O}_{n\times s\times p}$ is the tensor of  size $(n\times s\times p)$ having all its entries equal to zero.
\end{proposition}

\medskip
\noindent In the sequel, we briefly present  the tensor discrete cosine  GMRES algorithm to solve the problem (\ref{tr2}). %It could be considered as generalization of the  global GMERS algorithm \cite{jbilou1}. 
Let $\mathscr{  {X}}_{0}\in \mathbb{R}^{n\times s\times p}$ be an arbitrary initial guess with   the corresponding  residual
$\mathscr{R}_0=\mathscr{C}-\mathscr{A}\star_c \mathscr{X}_0$.    The aim of tensor  cosine GMRES method is to find and approximate solution  $\mathscr{X}_{m}$ approximating the exact solution $\mathscr{X}^*$   such that 

\noindent   \begin{equation}\label{xm}\mathscr{X}_{m}=\mathscr{X}_{0}+\mathbb{V}_{m}\circledast y, \end{equation}   
where $y=y_{m,\lambda_m} \in \mathbb{R}^m $ solves the projected regularized minimization problem 

\begin{eqnarray}\label{leastsq}
 	y_{m,\lambda_m}&=& \arg \min_{ y \in \mathbb{R}^{m }}\left(\Vert  \beta e_1-\widetilde{ { H}}_{m}  y   \|_ 2^2 +\lambda_m^2 \| y \|_2^2
 	\right),\\
 	&=&\arg \min_{ y \in \mathbb{R}^{m}}\left  \Vert \left ( \begin{array}{ll}
 	\widetilde H_m\\
 	\lambda_m I_m
 	\end{array}\right ) y - \left ( \begin{array}{ll}
 	\beta e_1\label{min}\\
 	0
 	\end{array}\right ) 
 	\right \Vert_2^2,
 	\end{eqnarray}
 	where $\beta=\|\mathscr{R}_0\| $ and $e_{1}$  the first  canonical basis vector in $\mathbb{R}^{m+1}$. 
 The minimizer $y_{m,\lambda_m}$ can also be computed as the solution of the following normal equations associated with  \eqref{min}
  \begin{equation}{\label{tikho2}}
  \widetilde H_{m,\lambda_m} y=\widetilde H_m^T\beta e_1, \quad \widetilde H_{m,\lambda_m}= (\widetilde H_m^T \widetilde C_m+ \lambda_m^2 I_m).
  \end{equation}
 		Note that since the Tikhonov problem \eqref{tikho2} is now a matrix one with small dimension as $m$ is generally small,  $\lambda_m$, can thereby be inexpensively computed by some techniques such as the GCV method \cite{golubwahba} or the L-curve criterion \cite{reichel2,reichel1,hansen1}. In this paper we consider the
 			generalized cross-validation (GCV) method
 	 	to choosing the regularization parameter 
  	\cite{golubwahba,wahbagolub}. 
  	We take advantage of the SVD decomposition of the low dimensional matrix $\widetilde H_m$ to obtain a more simple and computable expression of $GCV(\lambda_m)$. Consider the SVD decomposition  $\widetilde C_k=U\Sigma V^T$. Then, the GCV  function can be expressed as (see \cite{wahbagolub})
 	\begin{equation}
  	\label{gcv2}
  	GCV(\lambda_m)=\frac{\displaystyle
  		\sum_{i=1}^m\left(\frac{\widetilde
  			g_i}{\sigma_i^2+\lambda_m^2}\right)^2}{\displaystyle\left(\sum_{i=1}^m
  		\frac{1}{\sigma_i^2+\lambda_m^2}\right)^2},
  	\end{equation}
  	where $\sigma_i$ is the $i$th singular value of the matrix
 	 	$\widetilde H_m$ and $\widetilde g= \beta_1 U^T e_1$. The restarted tensor discrete cosine GMRES algorithm   is summarized as follows:
 		\begin{algorithm}[h!]
 		\caption{  Restarted tensor discrete cosine  GMRES (DC-GMRES(m)) method with  Tikhonov regularization  }\label{DC-GMRES}
 		\begin{enumerate}
 			\item 	{\bf Input.} $\mathscr{A}\in \mathbb{R}^{n\times n \times p}$, $\mathscr{C},\mathscr{X}_{0}\in \mathbb{R}^{n\times s \times p}$, an integer $m$ for restarting, a maximum number of iterations  $\text{Iter}_{\text{max}} $ and a tolerance $tol>0$ .
 			\item 	{\bf Output.} $  \mathscr{X}_{m}\in \mathbb{R}^{n\times s\times p}$ approximate  solution of the system  (\ref{eq1}).
 			%\item $k=1,\ldots,\text{Iter}_{\text{max}}$
 			\item  Compute $\mathscr{R}_{0}=\mathscr{C}-\mathscr{A}\star\mathscr{X}_{0} $, set $k=0$ .
 						%\begin{enumerate}
 				%\item  Compute $\mathscr{R}_{0}=\mathscr{C}-\mathscr{A}\star\mathscr{X}_{0} $.
 				\item  Apply Algorithm \ref{TGA}  to the pair $ (\mathscr{A},\mathscr{R}_0) $ to compute  $\mathbb{V}_{m}$ and  ${ \widetilde{H}}_m$ .
 				%			
 				%			\item Calculate the T-QR decomposition of  $\mathscr{ \widetilde{H}}_m$ using algorithm (\ref{TQRsimple}) . 
 				%			\item  Calculate $(\mathscr{R}_m)$ and $\mathscr{G}_m$
 				%			using the relations (\ref{grm}),(\ref{grm1}).
 				\item  Determine $\lambda_{m}$  as the parameter minimizing the GCV function  given by (\ref{gcv2})
 				\item  Compute the regularized solution $y_{m,\lambda_m}$ of the problem \eqref{min}. 
 				%	$(\mathscr{R}_m) \star \mathscr{Y}_m=   \mathscr{G}_m$ . 
 				\item Compute the approximate solution  $\mathscr{X}_{m}=\mathscr{X}_{0}+ \mathbb{V}_{m} \circledast y_{m,\lambda_m}  
 				$
 		%	\end{enumerate}
 			\item If $\|\mathscr{R}_{m}\|_F<tol$ or $k> itermax$, stop, \\
 			else			  Set    $\mathscr{X}_{0}=\mathscr{X}_{m}, k=k+1$ and go to Step 4.
 		%	\item End
 		\end{enumerate}
 	\end{algorithm}

	\subsection{The discrete cosine Golub-Kahan method}

	We consider the tensor least squares problem 
		\begin{equation}\label{ls1}
	\underset{\mathscr{X}}{\text{min}}\{\|\mathscr{A} \star_c \mathscr{X} - \mathscr{C}  \|_F^2,
	\end{equation}
where  $\mathscr{A} \in \mathbb{R}^{n\times \ell\times p}$  and   $\mathscr{C}  \in \mathbb{R}^{n\times  s \times p}$.  Instead of using the tensor cosine Arnoldi, we can use a dicrete cosine version of the tensor  Lanczos process to generate a new basis that can be used for the projection. We will use the tensor   Golub Kahan algorithm related to the c-product. 
%We notice here that we already defined in \cite{Elguide} another version of the tensor Golub Kahan algorithm by using the $m$-mode or the Einstein products with applications to color image restoration.\\
  defined as follows.
	
	\begin{algorithm}
		\caption{ The Tensor  discrete cosine Golub Kahan algorithm}\label{TG-GK}
		\begin{enumerate}
			\item {\bf Input.} The tensors $\mathscr {A}$,   $\mathscr {C}$ and an integer $m$.
			\item 	Set $\beta_1= \Vert \mathscr{C}\Vert_F$, $\alpha_1= \Vert \mathscr {A}^T\star_c\mathscr {U}_1  \Vert_F$, $\mathscr {U}_1=\mathscr {C}/\beta_1$ and 
			$\mathscr {V}_1=(\mathscr {A}^T\star_c\mathscr {U}_1) /\alpha_1$.
			\item for $j=1,\ldots,m$
			\begin{enumerate}
				\item $ \widetilde{\mathscr {U}}_{}= \mathscr {A} \star_c \mathscr {V}_{j} -\alpha_{j}\mathscr {U}_{j}$
				\item $\beta_{j+1}=\Vert \widetilde{ \mathscr {U}}_{}\Vert_F$ %if $\beta_j=0$ stop, else
				\item $\mathscr {U}_{j+1}=\widetilde {\mathscr {U}}/\beta_{j+1}$
				\item $\widetilde {\mathscr {V}}=\mathscr {A}^T \star_c \mathscr {U}_{j+1}-\beta_{j+1} \mathscr{V}_{j}$
				\item $\alpha_{j+1}=\Vert \widetilde {\mathscr {V}} \Vert_F$
				%\item if $\alpha_j=0$ stop, else
				\item $\mathscr {V}_{j+1}=\widetilde {\mathscr {V}}/\alpha_{j+1}$.
			\end{enumerate}
		\end{enumerate}
	\end{algorithm}

	\medskip 
\noindent Let $\widetilde{C}_m$ be the upper bidiagonal $((m+1) \times m  )$ matrix 
	$$ \widetilde{ {   {C}}}_m=\left[ \begin{array}{*{20}{c}}
	{{{\alpha}_1  }}&{{ }}& &   \\
	{\beta}_{2}&{{{\alpha}_2}}&\ddots& \\
	&\ddots&\ddots& \\
	&   & {\beta}_{m}  &  {\alpha}_{m}\\
	&    &       &    {\beta}_{m+1}
	\end{array}  \right] 
	$$
	and let $ {{{C}}}_m$ be the $(m \times m )$ matrix obtain   by deleting  the last row of  $\widetilde{{   {C}}}_m$. We denote by  $C_{.,j}$  the $j$-th column of the matrix  $C_m$. Let  $\mathbb{V}_{m}  $ and $\mathscr{A}\star_c\mathbb{V}_{m}  $ be the $(\ell\times (sm)\times p)$ and   $(n\times (sm)\times p)$ tensors with frontal slices $\mathscr{V}_{1},\ldots,\mathscr{V}_{m}$ and  $\mathscr{A}\star_c\mathscr{V}_{1},\ldots,\mathscr{A}\star_c\mathscr{V}_{m}$, respectively, and let  $\mathbb{U}_{m}  $ and $\mathscr{A}^T\star_c\mathbb{U}_{m}  $ be the $(n\times (sm)\times p)$ and $(\ell\times (sm)\times p)$  tensors with frontal slices $\mathscr{U}_{1},\ldots,\mathscr{U}_{m}$ and  $\mathscr{A}^T\star_c\mathscr{U}_{1},\ldots,\mathscr{A}^T\star_c\mathscr{U}_{m}$, respectively. We set  
	\begin{align}
	\label{ev12}
	\mathbb{U}_{m}:&=\left[  \mathscr{U}_{1},\ldots,\mathscr{U}_{m}\right], \;\;\; {\rm and}\;\;\; \mathscr{A}\star_c\mathbb{V}_{m}:=[\mathscr{A}\star_c\mathscr{V}_{1},\ldots,\mathscr{A}\star_c\mathscr{V}_{m}],\\
	\mathbb{V}_{m}:&=\left[  \mathscr{V}_{1},\ldots,\mathscr{V}_{m}\right], \;\;\; {\rm and} \;\;\; \mathscr{A}^T\star_c\mathbb{U}_{m}:=[\mathscr{A}^T\star_c\mathscr{U}_{1},\ldots,\mathscr{A}^T\star_c\mathscr{U}_{m}].
	\end{align}

	\begin{proposition}\label{proptggkb} 
		The tensors produced by the tensor cosine  Golub-Kahan algorithm satisfy the following relations
		\begin{eqnarray} \label{equa20}
	 			\mathscr{A} \star_c \mathbb{V}_m& = &\mathbb{U}_{m+1} \circledast {\widetilde { {   {C}}}}_m   ,    \\
			%	\begin{eqnarray} \label{equa30}
			\mathscr{A}^{T}\star_c\mathbb{U}_{m}& = &\mathbb{V}_m  \circledast {\widetilde { {   {C}}}}_m^T , \\
			\mathbb {U}_{m+1}  \circledast (\beta_1 e_1)&=&\mathscr {C},\\
			\Vert  \mathbb {U}_{m+1}  \circledast z \Vert_F & = & \Vert z \Vert_2,
			\end{eqnarray}
			where $e_1=(1,0,\ldots,0)^T\in \mathbb{R}^{m+1}$ and $z$ is a vector of  $\mathbb{R}^{m+1}$.
		\end{proposition} 
		
\medskip 

	\noindent To solve the least  squares problem  \eqref{ls1}, we consider approximations 	defined as
 \begin{equation}\label{xmlsqr}
 \mathscr{X}_{m}= \mathbb{V}_{m}\circledast y_m,  
 \end{equation} 
satisfying the minimization property of the corresponding residual. As we explained earlier, the problems that we are concerned with are ill-posed problems and then regularization techniques are highly recommended in those cases. But as the problem is very large, we apply  the regularization process to the projected problem derived from the minimization of the residual. This leads to a low dimensional Tikhonov formulation and then we seek for $y={y}_m \in \mathbb{R}^m $  that solves the low dimensional linear system of equations 
	\begin{equation}\label{normeq2}
	(\widetilde{C}_m^T\widetilde{C}_m+\lambda_m^2  I_m)y=\alpha_1\widetilde{C}_m^Te_1,\qquad\alpha_1=\|\mathscr{C}\|_F,
	\end{equation}
which is also equivalent to solving the least-squares problem
	\begin{equation}\label{leastsq0}
	\min_{y\in\mathbb{R}^m} \begin{Vmatrix}
	\begin{bmatrix}
	\lambda_m\widetilde{C}_m\\
	I_m
	\end{bmatrix}
	y-\alpha_1\lambda_me_1 \end{Vmatrix}_2\end{equation}
	
	\noindent  The regularized parameter $\lambda_m$ is computed by using the GCV function given by \eqref{gcv2}. The following algorithm summarizes the main steps of the  described method. 
%	\newpage
	
	\begin{algorithm}[!h]
		\caption{The  Tensor Discrete Cosine  Golub-Kahan (DC-GK) method}\label{TG-GKB}
		\begin{enumerate}
			\item {\bf Input.} The tensors $\mathscr{A}$,  $\mathscr {C}$.
			%\item {\bf Output.} cosine-GGKB steps $k$, $\lambda_{k}$ and $X_{k,\lambda_k}$.
			\item Determine the orthonormal bases $\mathbb{U}_{m+1}$ and $\mathbb{V}_{m}$ of tensors, and the bidiagonal $C_m$ and $\widetilde{C}_m$
			matrices  with Algorithm \ref{TG-GK}.
			\item Determine $\lambda_{m}$ using GCV function.
			\item  Determine $y_{m,\lambda_m}$ by solving  (\ref{leastsq0}) and then compute $X_{m,\lambda_m}$ by (\ref{xmlsqr}).
		\end{enumerate}
			\end{algorithm}
		In the next section, we derive a direct computation of the approximate Golub-Kahan solution by using a discrete cosine LSQR algorithm.
		
		\subsection{The  discrete cosine-LQSR method}
		\medskip
		\noindent In  this  section, we develop the tensor version of 
		the  well know LSQR algorithm introduced in \cite{Golub0} based on c-product formalism. 
		Let $\mathscr{A} \in \mathbb{R}^{n\times \ell \times p}$  be a tensor and  let  %$\mathscr{B}  \in \mathbb{R}^{\ell\times  s \times p}$ and  
		$\mathscr{C}  \in \mathbb{R}^{n\times  s \times p}$   a starting tensor.  
		
		\noindent The  purpose of  the  tensor  cosine  LSQR method is   to find, at some step $k$,   an approximation $ {\mathscr{X}}_{k}$ of the solution $ {\mathscr{X}}^*$  of  the  problem  (\ref{ls1}) such  that 
			\begin{equation}\label{sollsqr}
			\mathscr{X}_k= \mathbb{V}_k\circledast y_k,
			\end{equation} 
			where $y_k \in \mathbb{R}^k$. The associated    residual  is  given  by 
			\begin{equation}
			\mathscr{R}_k=  \mathscr {C}- \mathscr{A}\star_c\mathscr{X}_k =\beta_{1}\mathscr{U}_1-\mathbb{U}_{k+1}\circledast \widetilde { {   {C}}}_k \circledast   y_k=\mathbb{U}_{k+1} \circledast (\beta_{1}e_1-  \widetilde { {   {C}}}_k \circledast   y_k)
			\end{equation}
		and  using Proposition \ref{proptggkb}, we  get
			$$ \|\mathbb{U}_{k+1} \circledast (\beta_{1}e_1-  \widetilde { {   {C}}}_k \circledast   y_k)\|_F =\|\beta_{1}e_1-  \widetilde { {   {C}}}_k \circledast   y_k\|_2.$$ 
			\noindent This minimization problem is accomplished by using the QR  decomposition, where  a   unitary matrix $Q_k$ is determined so that  
			$$ Q_k  [    \widetilde { {   {C}}}_k \;\;\;  \beta_{1}e_1 ]=\left[ \begin{array}{*{20}{c}}
			{{ R_k  }}&{{ f_k}}   \\
			%{\beta}_{2}&{{{\alpha}_2}}&\ddots& \\
			%&\ddots&\ddots& \\
			%&   & {\beta}_{m}  &  {\alpha}_{m}\\
			0  &  \bar{\phi}_{k+1}
			\end{array}  \right] =\left[ \begin{array}{*{20}{c}}
			{{{\rho}_1  }}&{{\theta _2 }}&  & &  & \phi_{1} \\
			&{\rho}_{2}&{{{\theta}_3}}& & & \vdots\\
			&   &  \ddots   & \ddots    &    &\vdots \\
			& & & {\rho}_{k-1} & {\theta}_{k} & {\phi}_{k-1}\\
			&   & & & {\rho}_{k}  &  {\phi}_{k}\\
			&  &  &  &     &    \bar{\phi}_{k+1}
			\end{array}  \right] ,$$
			where $\rho_{i}, \theta{i}$ are  scalars. The  matrix $Q_k$  is  a  product  of  plane rotations designed  to  eliminate  the  sub-diagonals of  $\widetilde { {   {C}}}_k$.  This gives the following simple recurrence relation  : 
			$$  \left[ \begin{array}{*{20}{c}}
			{{ c_k  }}&{{ s_k}}   \\
			%{\beta}_{2}&{{{\alpha}_2}}&\ddots& \\
			%&\ddots&\ddots& \\
			%&   & {\beta}_{m}  &  {\alpha}_{m}\\
			- s_k  &  c_k
			\end{array}  \right]  \left[ \begin{array}{*{20}{c}}
			{{ \bar{\rho}_k  }}& 0 &{ \bar{ \phi}_{k}}   \\
			%{\beta}_{2}&{{{\alpha}_2}}&\ddots& \\
			%&\ddots&\ddots& \\
			%&   & {\beta}_{m}  &  {\alpha}_{m}\\
			\beta_{k+1}   &  \alpha_{k+1} & 0
			\end{array}  \right]=  \left[ \begin{array}{*{20}{c}}
			{{  {\rho}_k  }}&  \theta_{k+1} &{  { \phi_k}}   \\
			%{\beta}_{2}&{{{\alpha}_2}}&\ddots& \\
			%&\ddots&\ddots& \\
			%&   & {\beta}_{m}  &  {\alpha}_{m}\\
			0   &  \bar{\rho}_{k+1} & \bar{ \phi}_{k+1}
			\end{array}  \right], $$ 
			where $ \bar{\rho}_{1}=\alpha_{1}$ and $\bar{ \phi}_{1}=\beta_1$  and  the  scalars  $s_k, c_k $   are  the nontrivial element of $Q_{k+1,k}$ the  $k$-th plane   rotation. We  get 
			$$R_k y_k=f_k $$ and the  approximate  solution  is  given  by : $$\mathscr{X}_k= (\mathbb{V}_k \circledast R_k^{-1})\circledast f_k.$$ 
			Let $\mathbb{V}_k \circledast R_k^{-1}=\mathbb{P}_k=[\mathscr{P}_1 \ldots \mathscr{P}_k]$,  then  we  have $$ \mathscr{X}_k= \mathbb{P}_k   \circledast f_k$$
			\noindent  Notice  that the  tensor    $\mathscr{P}_k$ can  be  computed from    $\mathscr{P}_{k-1} $ and $\mathscr{V}_{k} $ as follows : $$\mathscr{P}_{k}   =(\mathscr{V}_{k}- \theta_{k}\mathscr{P}_{k-1}    ) \rho_{k}^{-1}.  $$  
			We also   have $f_k = \left[ \begin{array}{*{20}{c}}
			{{ f_{k-1}  }}   \\
			%{\beta}_{2}&{{{\alpha}_2}}&\ddots& \\
			%&\ddots&\ddots& \\
			%&   & {\beta}_{m}  &  {\alpha}_{m}\\
			\phi_k   
			\end{array}  \right] $ in  which $ \phi_k= c_k \bar{\phi}_k $. Finally,  $\mathscr{X}_k$ can  be  computed  as  follows $$\mathscr{X}_k= \mathscr{X}_{ k-1  }  + \phi_{k}\mathscr{P}_{k}.$$ Furthermore, we have  $$ \|\mathscr{R}_{k}  \|_F=|\bar{\phi}_{k+1}  |. $$
			\newpage 
			\noindent     The next algorithm which is named Discrete Cosine LSQR (DC-LSQR)  algorithm, describes the whole process.
		\begin{algorithm}[h!]
			\caption{ The Discrete Cosine LSQR (DC-LSQR)  algorithm    }\label{lsqr2}
			\begin{enumerate}
				\item {\bf Input.} The tensors $\mathscr {A}$,   $\mathscr {C}$,  $\mathscr {X}_0= \mathscr {0}$, $itermax$, the maximum number of allowed iterations and a tolerance $tol$ for the stopping criterion.
				\item 	Set $\beta_1= \Vert \mathscr{C}\Vert_F$, $\alpha_1= \Vert \mathscr {A}^T\star_c\mathscr {U}_1  \Vert_F$, $\mathscr {U}_1=\mathscr {C}/\beta_1$ and 
				$\mathscr {V}_1=(\mathscr {A}^T\star_c\mathscr {U}_1) /\alpha_1$, 	$\mathscr{W}_1=\mathscr{V}_1$ $, \bar{\rho}_{1}=\alpha_{1}$ and $\bar{ \phi}_{1}=\beta_1$ 
				\item for $j=1,\ldots, itermax$
				\begin{enumerate}
					\item $  {\mathscr {W}}_{j}= \mathscr {A} \star_c \mathscr {V}_{j} -\alpha_{j}\mathscr {U}_{j}$, 
				 $\beta_{j+1}=\Vert  { \mathscr {W}}_{j}\Vert_F$ and 
				 $\mathscr {U}_{j+1}=  {\mathscr {W}_j}/\beta_{j+1}$.
					\item $\widetilde {\mathscr {V}}=\mathscr {A}^T \star_c \mathscr {U}_{j+1}-\beta_{j+1} \mathscr{V}_{j}$, 
					$\alpha_{j+1}=\Vert \widetilde {\mathscr {V}} \Vert_F$ and  $\mathscr {V}_{j+1}=\widetilde {\mathscr {V}}/\alpha_{j+1}$.
					\item $\rho_{j}=(\bar{\rho}^2_j+\beta_{j+1}^2)^{\frac{1}{2}}$,  $ c_j=  \frac{\bar{\rho}_j}{\rho_j}$ and 
 $ s_j= \frac{\beta_{j+1}}{\rho_j}$.
					
					\item $ \theta_{j+1}=  {s_{j}}{\alpha_{j+1}}$ and  $ \bar{\rho}_{j+1}=  {c_{j}}{\alpha_{j+1}}$.
					\item $  {\phi}_{j}=  {c_{j}}\bar{\phi_{j}}$ and   $  \bar{\phi}_{j+1}=  -{s_{j}}\bar{\phi_{j}}$.
					\item $ \mathscr{X}_i=\mathscr{X}_{i-1}+ \frac{\phi_{j}}{\rho_{j}}   {\mathscr {W}}_j$;  $  {\mathscr {W}}_{j+1}=\mathscr{V}_{j+1} - \frac{\theta_{j+1}}{\rho_{j}}  {\mathscr {W}}_j$.
					\item If $|\bar{\phi}_{j+1}  | < tol$   stop. 
				\end{enumerate}
			\end{enumerate}
		\end{algorithm}
	
	\medskip 
	\noindent 	For ill posed problems, as it is the case for image or video restorations, we could have situations where the residual norm is small enough but the error norm is still large. As it is observed for those problems, the residual and the error norms could decrease in DC-LSQR till some iteration $k$ and then the norm of the error becomes to increase. One possibility to overcome these situations is to stop the iterations at some optimal $k_{opt}$.  The  L-curve criterion
		\cite{reichel2,hansen1} could be usefull to determine such optimal index $k_{opt}$. The method suggests to plot the curve $(\Vert \mathscr {R}_k \Vert,\Vert \mathscr {X}_k \Vert )$. Intuitively,
		the best regularization parameter should lie on the corner of the
		L-curve  corresponding to the point on the curve with maximum curvature.

	\section{Numerical results}
	In this section, we give some  numerical tests on the methods described in this paper. We compared the performances of the   tensor discrete cosine GMRES describes in Algorithm \ref{DC-GMRES},  the tensor discrete cosine Golub-Kahan (DC-GK) algorithm given by  Algorithm \ref{TG-GKB} and the tensor dicrete cosine LSQR (DC-LSQR) described in Algorithm  \ref{lsqr2} ,   
	when applied to the restoration of blurred and noisy color images.  All computations were carried out using the Matlab environment on an Intel(R) Core(TM) i7-8550U CPU @ 1.80GHz (8 CPUs) computer with 12 GB of
	RAM. The computations were done with approximately 15 decimal digits of relative
	accuracy. 
	Let $\widehat{X}^{(1)}$, $\widehat{X}^{(2)}$, and $\widehat{X}^{(3)}$ be the $n\times n$ matrices that constitute the three channels of the original error-free color image $\widehat{\mathscr{X}}$, and $\widehat{C}^{(1)}$, $\widehat{C}^{(2)}$, and $\widehat{C}^{(3)}$ the $n\times n$ matrices associated with error-free blurred color image $\widehat{\mathscr{C}}$. 
	We consider that both cross-channel and within-channel blurring  take place in the blurring process of the original image. The $\tt{vec}$ operator 
	transforms a matrix  to a vector  by stacking the columns of the matrix from
	left to right. The full blurring model can be descried as follows 
	\begin{equation}\label{linmodel}
	\left(\mathbf{A}_{\text {color }} \otimes \mathbf{A^{(1)}}\otimes\mathbf{A^{(2)}}\right) \widehat{\mathbf{x}}=\widehat{\mathbf{c}},
	\end{equation}
	where, 
	$$ \widehat{\mathbf{c}}=\left[\begin{array}{c}
	\tt{vec}\left(\widehat{\mathbf{C}}^{(1)}\right) \\
	\tt{vec}\left(\widehat{\mathbf{C}}^{(2)}\right) \\
	\tt{vec}\left(\widehat{\mathbf{C}}^{(3)}\right)
	\end{array}\right], \quad \widehat{\mathbf{x}}=\left[\begin{array}{c}
	\tt{vec}\left(\widehat{\mathbf{X}}^{(1)}\right) \\
	\tt{vec}\left(\widehat{\mathbf{X}}^{(2)}\right) \\
	\tt{vec}\left(\widehat{\mathbf{X}}^{(3)}\right)
	\end{array}\right]\; and \; 	\mathbf{A}_{\mathrm{color}}=\left[\begin{array}{ccc}
	a_{\text{rr}} & a_{\text{rg}}& a_{\text{rb}} \\
	a_{\text{gr}} & a_{\text{gg}} & a_{\text{gb}} \\
	a_{\text{br}} & a_{\text{bg}} & a_{\text{bb}}
	\end{array}\right],$$
	where $\mathbf{A}_{\text {color }}$ is the $3\times3$ matrix obtained from \cite{HNO}, that models the  cross-channel blurring, in which  each row sums is  one. We consider the special case where $a_{\text{rr}}=a_{\text{gg}}=a_{\text{bb}}$, $a_{\text{gr}}=a_{\text{rg}}$, $a_{\text{br}}=a_{\text{rb}}$, and $a_{\text{bg}}=a_{\text{gb}}$, which gives rise to a cross-channel  circular mixing.  $\mathbf{A^{(1)}}\in\mathbb{R}^{n\times n}$ and $\mathbf{A^{(2)}}\in\mathbb{R}^{n\times n}$ define within-channel blurring and they model
	the horizontal within blurring and the vertical  within blurring matrices, respectively; for more details see \cite{HNO} where the  notation $\otimes$ stands for  the Kronecker product of  matrices. By exploiting  the circulant structure  of the cross-channel blurring matrix $\mathbf{A}_{\text {color }}$, it can be easily shown that (\ref{linmodel}) can be written in the following tensor form
	\begin{equation}\label{tensor_form}
	\mathscr {A}\star_c \widehat{\mathscr{X}}\star_c \mathscr{B}= \widehat{\mathscr{C}},
	\end{equation}
	where $\mathscr {A}$ is a 3-way tensor such that  $\mathscr {A}(:,:,1)=\alpha \mathbf{A^{(2)}}$, $\mathscr {A}(:,:,2)=\beta \mathbf{A^{(2)}}$ and $\mathscr {A}(:,:,3)=\gamma \mathbf{A^{(2)}}$ and  $\mathscr {B}$ is a 3-way tensor with $\mathscr {B}(:,:,1)=(\mathbf{A^{(1)}})^T$, $\mathscr {B}(:,:,2)=0$ and $\mathscr {B}(:,:,3)=0$. To test the performance of algorithms, the within blurring matrices $A^{(i)}$  have the following
	entries 
	$$a_{k \ell}=\left\{\begin{array}{ll}
	\frac{1}{\sigma \sqrt{2 \pi}} \exp \left(-\frac{(k-\ell)^{2}}{2 \sigma^{2}}\right), & |k-\ell| \leq r \\
	0, & \text { otherwise. }
	\end{array}\right.$$
	Note  that $\sigma$ controls the amount of smoothing, i.e. the larger the $\sigma$,  the more ill posed the problem.  We generated a blurred and noisy tensor image $\mathscr{C}=\widehat{\mathscr{C}}+\mathscr{N},$ where $\mathscr{N}$ is a noise tensor with normally distributed random entries with zero mean and with variance chosen to correspond to a specific noise level $\nu:=\|\mathscr{N}\|_F /\|\widehat{\mathscr{C}}\|_F.$
	To  compare the effectiveness of our solution methods, we evaluate 
	$$\text{Relative error}=\frac{\left\|\hat{ \mathscr X}-{\mathscr  X}_{\text{restored}}\right\|_{F}}{\|\widehat{ \mathscr X}\|_{F}}$$
	and the Signal-to-Noise
	Ratio (SNR) defined by
	\[\text{SNR}({ \mathscr X}_{\text{restored}})=10\text{log}_{10}\frac{\|\widehat{\mathscr X}-E(\widehat{\mathscr X})\|_F^2}{\|{\mathscr X}_{\text{restored}}-\widehat{\mathscr X}\|_F^2},\]
	where $E(\widehat{\mathscr X})$ denotes the mean gray-level of the uncontaminated image $\widehat{\mathscr{X}}$. 

\begin{figure}
	 	\begin{center}
	 			\includegraphics[width=1.5in]{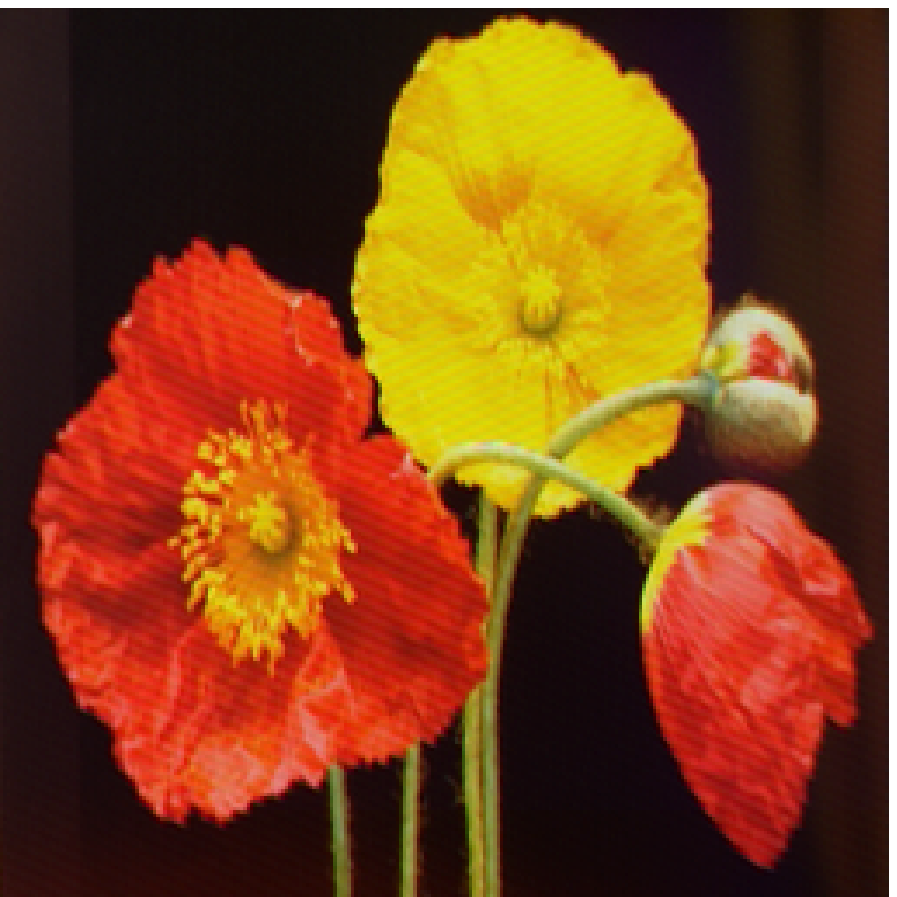}
	 			\includegraphics[width=1.5in]{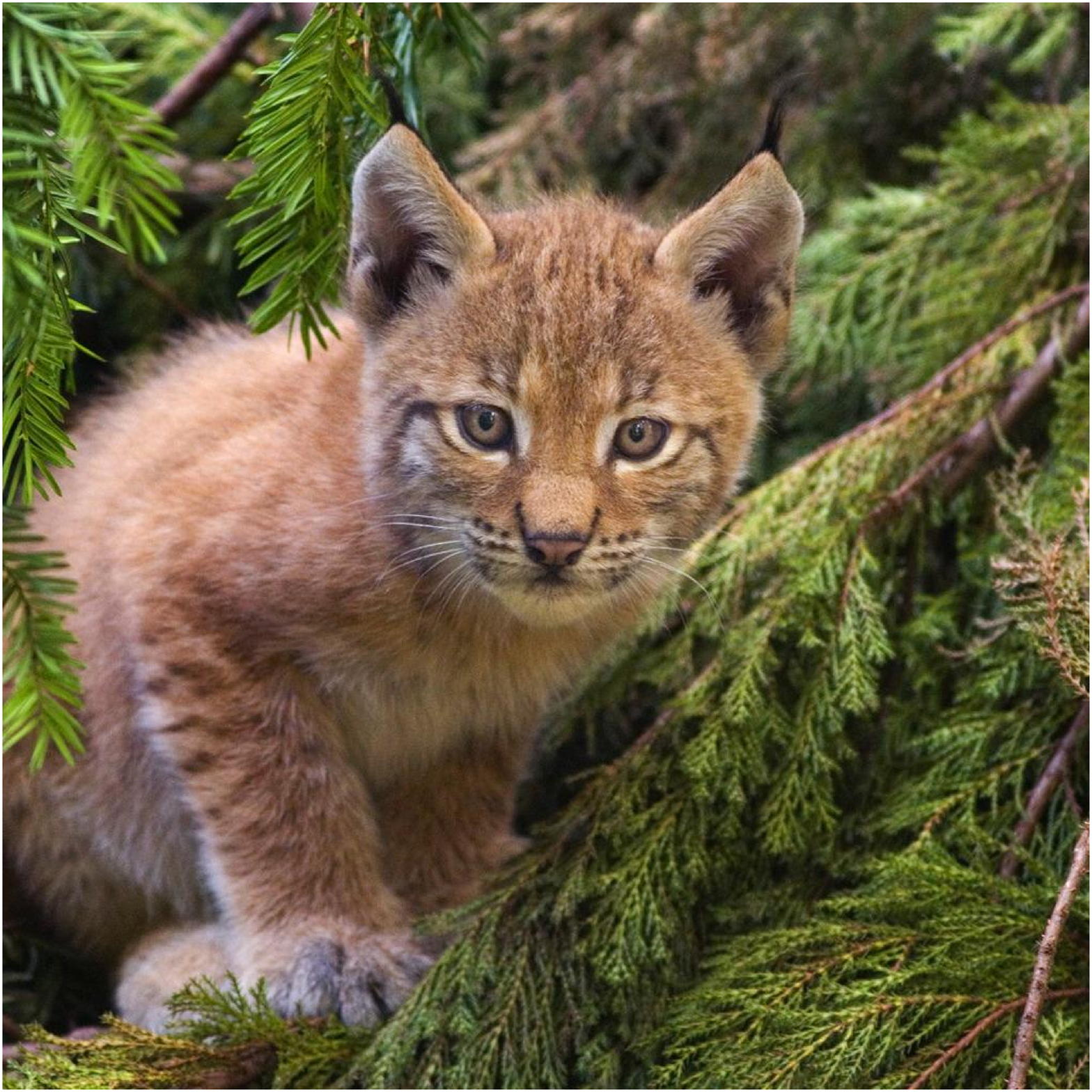}
	 		\caption{Original RGB images:  $\tt papav256$ (left) and $\tt cat1024$ (right).}\label{fig1}
	 	\end{center}
	 \end{figure}

	In our experiments,   we applied the three algorithms DC-GMRES(10), DC-GK and DC-LSQR  for the reconstruction of a cross-channel blurred color images that have been contaminated by both within and cross blur, and additive noise. The cross-channel blurring is determined by the matrix 
	$$\mathbf{A}_{\mathrm{color}}=\left[\begin{array}{ccc}
	0.8 & 0.10& 0.10 \\
	0.10 & 0.80 & 0.10 \\
	0.10& 0.10 & 0.80
	\end{array}\right].$$
	 We consider two $\mathrm{RGB}$ images, $\tt papav256$ ($\widehat{\mathscr X}\in\mathbb{R}^{256\times256\times3}$) and  $\tt cat1024$ ($\widehat{\mathscr X}\in\mathbb{R}^{1024\times 1024\times3}$). They are shown on  Figure \ref{fig1}. For the within-channel blurring,  we let $\sigma=4$ and $r=6$.  The associated blurred and noisy RGB images are obtained as $\mathscr{C}=\mathscr{A}\ast\widehat{\mathscr{X}}\ast\mathscr{B}+\mathscr{N}$.  Given the contaminated RGB image $\mathscr{C}$, we would like to recover an approximation of the original RGB image $\widehat{\mathscr X}$.   
	 
	  \subsection{Example 1} In the first experiment, we used the $\tt papav256$ color image of size $256  \times 256 \times 3$ with two different  noise levels $\nu=10^{-2}$ or $\nu = 10^{-3}$. In Table \ref{tab1} we reported the obtained SNR, the corresponding relative error norm and the required cpu-time for 
 DC-GMRES(10), DC-GK and DC-LSQR  with a noise level of $10^{-3}$.

\begin{figure}[h!]
	\begin{center}
		\includegraphics[width=5.5in]{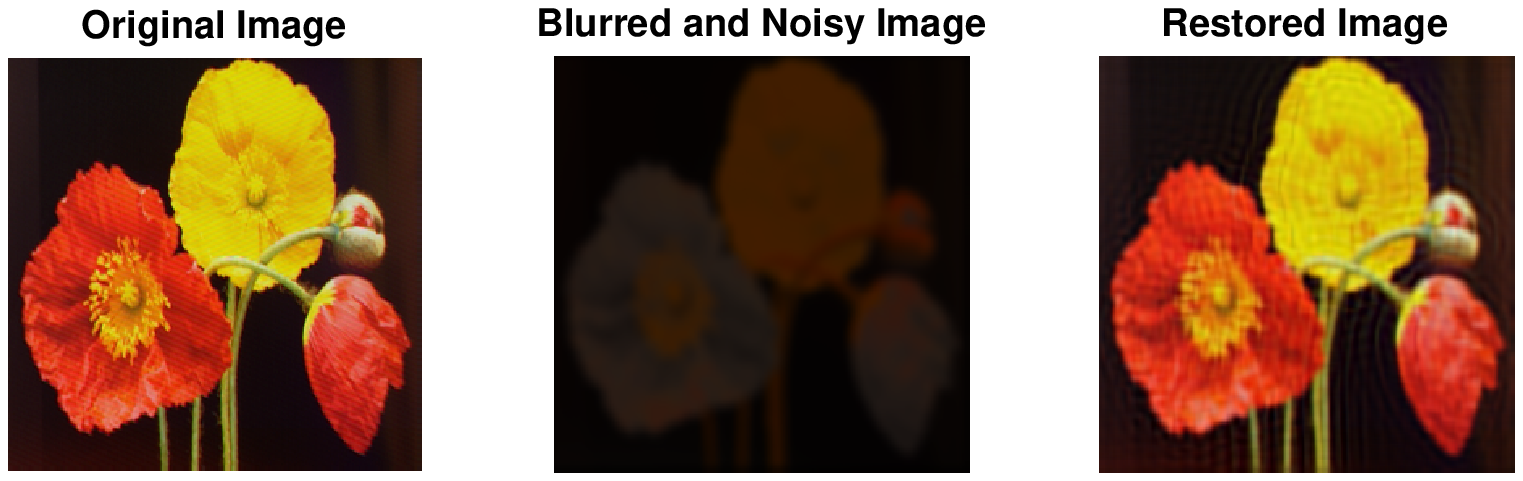}
		\caption{Test for Example 1, with DC-LSQR for $\tt papav256$, and  noise level $10^{-3}$. Original (left), noisy-blurred (center) and restored (right) .}\label{fig2}
	\end{center}
\end{figure}

\begin{table}[htbp]
	\begin{center}
			\caption{Results for Experiment 1 with $\tt papav256$. Noise level $10^{-3}$.}\label{tab1}
		\begin{tabular}{lcccc}
			\hline RGB images & Method & SNR & Relative error & cpu-time (seconds) \\
			\hline {$\tt papav256$}&DC-GMRES(10)&$19.25$& $8.5\times10^{-2}$&$5.21$ \\  
			&   DC-GK&$23.9$&$4.7\times10^{-2}$&$1.62$ \\
			&DC-LSQR&$21.8$&$6.8\times10^{-2}$&$1.23$\\
			\hline
		\end{tabular}
	\end{center}
\end{table}
 \noindent As can be seen from those results, the DC-LSQR requires lower cpu-time as compared to the other two methods. However, DC-GK returns the best SNR.  For this experiment, the optimal iteration number was $k_{opt}=14$. A maximum number of iterations was $itermax=10$ for DC-GMRES(10) and $mmax=15$ for DC-GK. As we mentioned earlier, DC-GMRES(10) and DC-GK were run with   the Tikhonov  regularization technique  (applied to the projected least squares problem) and we used GCV method for estimating the regularization parameters in each iteration of the processes. The obtained optimal parameters, at the final step were $\lambda_{1}=2.32\times 10^{-5}$ for DC-GMRES(10) and $\lambda_{1}=1.24\times 10^{-6}$ for DC-GK. 
Figure \ref{fig2} shows  the obtained blurred image and the restored one when using DC-LSQR method with noise level of $10^{-3}$. 

\noindent In Table \ref{tab2}, we reported the results obtained by DC-GMRES(10), DC-GK and DC-LSQR for the color image $\tt papav256$ with a noise level of $10^{-2}$. Here also, we used   $k_{opt}=15$ for DC-LSQR, $itermax=15$ for DC-GMRES(10) and $mmax=20$ for DC-GK. As can be seen, the DC-LSQR returns the best results when comparing the three methods. 

\begin{table}[htbp]
	\begin{center}
			\caption{Results for Example  1 with $\tt papav256$. Noise level $10^{-2}$.}\label{tab2}
		\begin{tabular}{lcccc}
			\hline RGB image & Method & SNR & Relative error & cpu-time (seconds) \\
			\hline {$\tt papav256$}&DC-GMRES(10)&$17.24$& $1.5\times10^{-2}$&$ 7.14$ \\  
			&   DC-GK&$20.4$&$7.2\times10^{-2}$&$1.92$ \\
			&DC-LSQR&$20.2$&$8.5\times10^{-2}$&$1.23$\\
			\hline
		\end{tabular}
	\end{center}
\end{table}

\subsection{Example 2} In the second example, we used the color image $\tt cat1024$ of dimension $1024 \times 1024 \times 3$. Here also, we compared the three methods using two noise levels $\nu=10^{-3}$ and $\nu=10^{-2}$. Table \ref{tab3} reports on the obtained results for the noise level $\nu=10^{-3}$ . For this experiment, the optimal iteration number for DC-LSQR was $15$, the maximum iteration number allowed to DC-GMRES(10) was $10$ and a maximum of $mmax=20$ iterations was for DC-BK. As can be seen from the obtained results, DC-LSQR returns the best results: for SNR and the total cpu-time. 

 \begin{table}[h!]
 	\caption{Results for Example 2 with  noise level $10^{-3}$.}\label{tab3}
 	\begin{center}
 		\begin{tabular}{lcccc}
 			\hline RGB image & Method & SNR & Relative error & cpu-time (seconds) \\
 			\hline 
 			 $\tt cat1024$ &  DC-GMRES(10)&$14.96$&$4.54\times10^{-2}$&$100.35$ \\
 			&DC-GK&$19.25$&$6.97\times10^{-2}$&$25.33$ \\
 			&DC-LSQR&$18.87$&$5.43\times10^{-2}$&$19.45$\\
 			\hline
 		\end{tabular}
 	\end{center}
 \end{table}

\begin{figure}[h!]
	\begin{center}
		\includegraphics[width=5.5in]{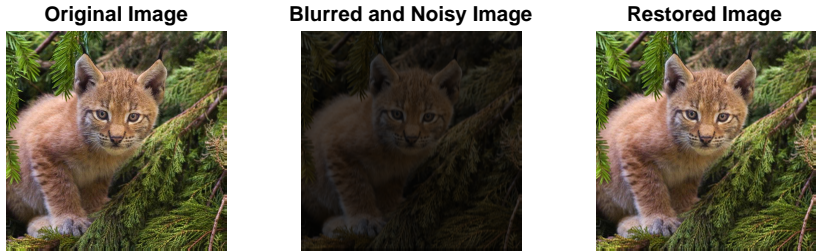}
\caption{Test for Example 2, with DC-LSQR for $\tt cat1024$, and  noise level $10^{-3}$. Original (left), noisy-blurred (center) and restored (right).}\label{fig3}
	\end{center}
\end{figure}
\noindent Figure \ref{fig3} shows  the obtained blurred image and the restored one when using DC-LSQR method with noise level of $10^{-3}$ for the image $\tt cat1024$. 
Table \ref{tab4} reports on the obtained results for the noise level $\nu=10^{-2}$.   For this experiment, the optimal iteration number for DC-LSQR was $20$, the maximum iteration number allowed to DC-GMRES(10) was $15$ and a maximum of $mmax=25$ iterations was for DC-BK. As can be seen from this table,  DC-LSQR returns the best results: for SNR and the total cpu-time. For the returned SNR, generally the two Golub Kahan based methods return similar results but the second formulation of the method which corresponds the DC-LSQR (Algorith \ref{lsqr2}) requires less cpu-time. 
\begin{table}[htbp]
	\begin{center}
			\caption{Results for Example 2 with  noise level $10^{-2}$.}\label{tab4}
		\begin{tabular}{lcccc}
			\hline RGB image & Method & SNR & Relative error & cpu-time (seconds) \\
			\hline 
			$\tt cat1024$ &  DC-GMRES(10)&$14.53$&$9.62\times10^{-2}$&$137.43$ \\
			&DC-GK&$15.87$&$8.06\times10^{-2}$&$30.22$ \\
			&DC-LSQR&$15.75$&$8.17\times10^{-2}$&$26.43$\\
			\hline
		\end{tabular}
	\end{center}
\end{table}

\section*{Conclusion} In this paper, we presented three discrete cosine Krylov-based methods, namely  tensor DC-GMRES, DC-GK and DC-LSQR. The second two methods use the discrete cosine Golub-Kahan bidiagonalisation algorithm  that we defined in this work. DC-GMRES and DC-GK are combined with the well known  Tikhonov regularization method that is applied, at each iteration for the two algorithms, to the obtained projected low dimensional ill-posed least squares minimisation problem. The reported numerical tests show that the methods are very fact and can be used as restoration techniques for color image restoration.

\end{document}